\documentclass[11pt,twoside,a4paper]{article}
\usepackage{amsmath}
\usepackage{amssymb}
\usepackage[latin1]{inputenc}
\usepackage[T1]{fontenc}   
\usepackage[francais]{babel}
\input{xy}
\xyoption{all}

\newcommand{\M}{{\cal M}}

\title{L'algèbre des invariants d'un groupe de Coxeter agissant sur un mutiple de
sa représentation standard}
\date{}
\author{L. Foissy \\
\\
{\small{\it Laboratoire de Math\'ematiques - UMR6056, Universit\'e de Reims}}\\
\small{{\it Moulin de la Housse - BP 1039 - 51687 REIMS Cedex 2, France}}\\
\small{e-mail : loic.foissy@univ-reims.fr}}

\newtheorem{prop}{\indent Proposition}
\newtheorem{defi}[prop]{\indent D\'efinition}
\newtheorem{lemme}[prop]{\indent Lemme}
\newtheorem{theo}[prop]{\indent Th\'eor\`eme}
\newtheorem{cor}[prop]{\indent Corollaire}

\begin{document}

\maketitle

RESUME : Soit $G$ un groupe de Coxeter de type $A_n$, $B_n$, $D_n$
ou $I_2(N)$, ou un groupe de réflexions complexes de type $G(de,e,n)$.
Soit $V$ sa représentation standard et soit $k$ un entier plus grand que $2$.
Alors $G$ agit sur $S(V)^{\otimes k}$.
Nous montrons que l'algèbre d'invariants $(S(V)^{\otimes k})^G$ est un 
$(S(V)^G)^{\otimes k}$-module libre de rang $|G|^{k-1}$ et que
$S(V)^{\otimes k}$ n'est pas un $(S(V)^{\otimes k})^G$-module libre.\\

ABSTRACT: Let $G$ be a Coxeter group of type $A_n$, $B_n$, $D_n$
or $I_2(N)$, or a complex reflection group of type $G(de,e,n)$.
Let $V$ be its standard representation and let $k$ be an integer greater than $2$.
Then $G$ acts on $S(V)^{\otimes k}$. 
We show that the algebra of invariants $(S(V)^{\otimes k})^G$ is a
free $(S(V)^G)^{\otimes k}$-module of rank $|G|^{k-1}$, and that
$S(V)^{\otimes k}$ is not a free $(S(V)^{\otimes k})^G$-module.\\

KEY-WORDS: Theory of invariants; Coxeter groups.\\

MSC CLASSES: 13A50, 17B10.

\tableofcontents

\section*{Introduction}

Soit $G$ un groupe de Coxeter et $V$ sa représentation standard.
Soit $A_1=S(V)$ l'algèbre symétrique de $V$ ; le groupe $G$ agit par automorphismes d'algèbre
sur $A_1$ et les résultats suivants sont bien connus (\cite{Bourbaki,Humphreys,Smith}) :
\begin{enumerate}
\item La sous-algèbre des éléments $G$-invariants $A_1^G$ est une algèbre de polynômes.
\item Le $A_1^G$-module $A_1$ est libre, de rang fini égal au cardinal de $G$.
\end{enumerate}

Soit $k$ un entier supérieur ou égal à $2$. 
On considère maintenant la représentation $V^{\oplus k}$ de $G$ et
son algèbre symétrique $A_k=S(V^{\oplus k})=A_1^{\otimes k}$, sur laquelle agit $G$. 
Pour $k=2$, cette situation a été étudiée dans le cas du groupe symétrique dans \cite{Haiman} ; 
en remarquant que $V$ et $V^*$ sont des représentations isomorphes, 
cette situation est étudiée dans \cite{Gordon} pour les groupes de Coxeter et
dans \cite{Alev} pour les groupes de Weyl de rang $2$.

Cette algèbre $A_k$ contient deux sous-algèbres particulières.
La première est la sous-algèbre $A_k^G$ des éléments $G$ invariants, la seconde
est la sous-algèbre des éléments $A_k^{G\times \ldots \times G}$ invariants, c'est-à-dire
$(A_1^G)^{\otimes k}$. De plus, $(A_1^G)^{\otimes k} \subseteq A_k^G \subseteq A_k$.
D'après les premiers résultats évoqués dans cette introduction,
$A_k$ est un $(A_1^G)^{\otimes k}$-module libre de rang fini, égal à $|G|^k$.
La question est de savoir si $A_k^G$ est un module libre sur $(A_1^G)^{\otimes k}$
et le cas échéant de calculer son rang. 
On peut également se demander si $A_k$ est un module libre sur $A_k^G$. 
Nous démontrons dans ce texte les résultats suivants :
\begin{enumerate}
\item Si $G$ est un groupe de Coxeter d'une des séries infinies $A_n$, $B_n$, $D_n$, $I_2(N)$
ou plus généralement un groupe de réflexions complexes de la série infinie $G(de,e,n)$, 
alors $A_k^G$ est un $(A_1^G)^{\otimes k}$-module libre, de rang $|G|^{k-1}$.
\item Avec les mêmes hypothèses sur $G$, $A_k$ n'est pas un $A_k^G$-module libre.
\end{enumerate}
Nous donnons également un exemple de groupe pour lequel $A_k^G$ n'est pas un $(A_1^G)^{\otimes k}$-module libre.
Ce groupe est le groupe diédral $I_2(6)=G_2$, mais la représentation choisie, de dimension 2, 
n'est pas la représentation standard. 

La preuve de ces résultats utilise une variante $\mathbb{N}^k$-graduée de la formule de Molien,
exposée dans la première section.
Les différents groupes considérés dans ce texte sont tous des sous-groupes
de produits en couronne de certains groupes cycliques, comme il est expliqué dans la deuxième section.
La section suivante détaille le calcul des séries de Poincaré-Hilbert des algèbres d'invariants 
et la dernière section explicite les différents exemples de cet article. \\

{\bf Notations.} 
\begin{enumerate}
\item Le corps de base est $\mathbb{C}$. 
\item Soit $A$ un anneau. Pour tout $q$ dans cet anneau, pour tout $i \in \mathbb{N}^*$, on pose
$[i]_q=1+\ldots+q^{i-1}$.
\end{enumerate}

\section{Formule de Molien $\mathbb{N}^k$-graduée}

Soit $G$ un groupe fini, agissant de manière homogène sur un espace vectoriel $\mathbb{N}^k$-gradué $A$. 
Les composantes homogènes de $A$ seront notées $A(i_1,\ldots,i_k)$.	
La série formelle de Poincaré-Hilbert de $A$ est :
$$R(h_1,\ldots,h_k)=\sum_{i_1,\ldots,i_k} dim_{\mathbb{C}}(A(i_1,\ldots,i_k))h_1^{i_1}\ldots h_k^{i_k}.$$

$G$ agit sur $A$ de manière homogène, donc le sous-espace $A^G$ des invariants sous l'action
de $G$ est un sous-espace gradué. On note $R^Gk(h_1,\ldots,h_k)$ sa série formelle de Poincaré-Hilbert.

\begin{defi}
\textnormal{
Soit $\sigma \in G$ et $(i_1,\ldots,i_k)\in \mathbb{N}^k$.
On pose $\chi_{i_1,\ldots,i_k}(\sigma)=Tr\left(\sigma_{\mid A(i_1,\ldots,i_k)}\right)$.
On pose également :
$$\chi_k(\sigma)=\sum_{i_1,\ldots,i_k}\chi_{i_1,\ldots,i_k}(\sigma)h_1^{i_1}\ldots h_k^{i_k}.$$
}
\end{defi}

On a la variante suivante de la formule de Molien (voir \cite{Smith})  :
\begin{prop}
\label{prop2}
La série formelle de l'espace $A^G$ des invariants de $A$ sous l'action de $G$ est :
$$R^G(h_1,\ldots,h_k)=\frac{1}{|G|} \sum_{\sigma\in G} \chi_k(\sigma).$$
\end{prop}

{\bf Preuve.} Il suffit de montrer que pour tout $(i_1,\ldots,i_k) \in \mathbb{N}^k$,
$$dim_{\mathbb{C}}(A(i_1,\ldots,i_k))=\frac{1}{|G|} \sum_{\sigma\in G}\chi_{i_1,\ldots,i_k}(\sigma),$$
ce qui est classique. $\Box$\\

{\bf Exemple.} Soit $V$ une représentation de $G$ de dimension finie et soit $k\in \mathbb{N}^*$. 
Alors $G$ agit par automorphismes d'algèbre sur $A_k=S(V)^{\otimes k}=S(V^{\oplus k})$. 
Cette algèbre est $\mathbb{N}^k$-graduée en mettant les éléments de
la $i$-ème copie de $V$ homogènes de degré $(0,\ldots,0,1,0,\ldots,0)$,
où le coefficient $1$ est situé en $i$-ème position. 
La série formelle de Poincaré-Hilbert de $A_k$ est alors :
$$R_k(h_1,\ldots,h_k)=\sum_{i_1,\ldots,i_k} dim_{\mathbb{C}}(A_k(i_1,\ldots,i_k))h_1^{i_1}\ldots h_k^{i_k}
=\frac{1}{(1-h_1)^n} \ldots \frac{1}{(1-h_k)^n}.$$

\section{cadre et énoncé du théorème principal}

Soient $N\geq 2$ et $n \in \mathbb{N}^*$. Soit $H$ un sous-groupe de $(\mathbb{Z}/N\mathbb{Z})^n$ stable
sous l'action de  $S_n$ par permutation des coordonnées. On obtient alors un produit semi-direct 
$G=H \rtimes S_n$, sous-groupe du produit en couronne $(\mathbb{Z}/N\mathbb{Z})^n \rtimes S_n$.
Soit $V=Vect(x_1,\ldots,x_n)$ et soit $\xi$ une racine $N$-ième primitive de l'unité.
Le groupe $G$ agit sur $V$ de la manière suivante : pour tout $(\overline{k}_1,\ldots,\overline{k}_n) \in H$,
tout $\sigma \in S_n$,
$$(\overline{k}_1,\ldots,\overline{k}_n).x_j=\xi^{k_j} x_j, \hspace{1cm} \sigma.x_j=x_{\sigma(j)}.$$
Le but de la section \ref{sect3} du présent texte est de démontrer le résultat suivant :

\begin{prop}
\label{prop3}
Sous les hypothèses exposées précédemment,
$$\lim_{(h_1,\ldots,h_k) \longrightarrow (1,\ldots,1)} \frac{R^G_k(h_1,\ldots,h_k)}{R^G_1(h_1)\ldots R^G_1(h_k)}
=|G|^{k-1}.$$
\end{prop}

Par la suite, nous noterons :
$$Q_k(h_1,\ldots,h_k)=\frac{R^G_k(h_1,\ldots,h_k)}{R^G_1(h_1)\ldots R^G_1(h_k)}.$$

La proposition \ref{prop3} a le corollaire suivant :
\begin{theo}
\label{theo4}
Supposons que $(G,V)$ soit un groupe de Coxeter d'une des séries infinies $A_n$, $B_n$, $D_n$ ou $I_2(N)$,
ou plus généralement un groupe de réflexions complexes de la série infinie $G(de,e,n)$. 
Notons $A_k=S(V)^{\otimes k}=S(V^{\oplus k})$.
\begin{enumerate}
\item Pour tout $k\in \mathbb{N}^*$,
$A_k^G$ est un $(A_1^G)^{\otimes k}$-module libre de rang $|G|^{k-1}$. 
De plus, $Q_k(h_1,\ldots,h_k)$ est un polynôme.
\item Si $k\geq 2$, $A_k$ n'est pas un module libre sur $A_k^G$.
\end{enumerate}
\end{theo}

{\bf Preuve.} Lorsque $(G,V)$ est un groupe de réflexions complexes, 
on sait que $A_1^G$ est un anneau de polynômes 
et que $A_1$ est un $A_1^G$-module libre de rang fini (voir par exemple \cite{Bourbaki,Humphreys}). 
Par suite, $A_k=A_1^{\otimes k}$ est un $(A_1^G)^{\otimes k}$-module libre de type fini. Soit alors
$(P_1,\ldots,P_m)$ une $(A_1^G)^{\otimes k}$-base de $A_k$. \\

{\it Première étape.} On note $\M$ l'idéal d'augmentation de $(A_1^G)^{\otimes k}$.
Montrons que $\M A_k \cap A_k^G=\M A_k^G$. L'inclusion $\supseteq$ est immédiate. 
Soit $x \in \M A_k \cap A_k^G$.
Cet élément peut s'écrire $x=\sum_j m_j a_j$, avec pour tout $j$, $m_j \in \M$ et $a_j \in A_k$.
De plus, comme $x$ et les $m_j$ sont $G$-invariants :
\begin{eqnarray*}
x&=&\frac{1}{|G|} \sum_{g \in G} g.x\\
&=&\frac{1}{|G|} \sum_{g \in G} \sum_{j}  (g.m_j)(g.a_j)\\
&=&\frac{1}{|G|} \sum_{g \in G} \sum_{j}  m_j(g.a_j)\\
&=&\sum_{j}  m_j\underbrace{\left(\frac{1}{|G|}\sum_{g \in G} g.a_j\right)}_{\in A_k^G} \in \M A_k^G.
\end{eqnarray*}
par suite, on a une injection de $\mathbb{C}$-espaces vectoriels :
$$\frac{A_k^G}{\M A_k^G} \hookrightarrow  \frac{A_k}{\M A_k}.$$
Une $\mathbb{C}$-base de $ \frac{A_k}{\M A_k}$ est $(P_1+\M A_k, \ldots, P_m+\M A_k)$, 
donc $\frac{A_k}{\M A_k}$ et par suite $\frac{A_k^G}{\M A_k^G}$ sont de dimension finie.
Soit $(Q_1+\M A_k^G,\ldots,Q_n +\M A_k^G)$ une $\mathbb{C}$-base de $\frac{A_k^G}{\M A_k^G}$.
Notons que $n \leq m$ et qu'on peut choisir les $Q_i$ homogènes.\\

{\it Deuxième étape.} Montrons que les $Q_i$ engendrent $A_k^G$. 
Posons $B=_{(A_1^G)^{\otimes k}}\langle Q_1,\ldots,Q_n\rangle$ et $A=A_k^G/B$. Tout d'abord, $\M A=A$.
En effet, si $x \in A_k^G$, alors il existe $\lambda_1,\ldots,\lambda_n \in \mathbb{C}$, tels
que $x+\M A_k^G=\lambda_1 Q_1+\ldots +\lambda_n Q_n +\M A_k^G$. On peut donc écrire :
$$x=\lambda_1 Q_1+\ldots +\lambda_n Q_n +\sum_j m_j a_j,$$
où $a_j \in A_k^G$, $m_j \in \M$ pour tout $j$. En conséquence :
$$x+B=0+\sum_j m_j(a_j+B) \in \M A.$$ 
Par conséquence, pour tout $k \in \mathbb{N}^*$, $A=\M^kA$. \\

Supposons $A$ non nul. Il s'agit d'un $(A_1^G)^{\otimes k}$-module gradué ; choisissons $x \in A$, non nul,
homogène et notons $k$ son degré. Alors $x \in \M^{k+1}A$, donc peut s'écrire :
$$x=\sum_j m_j^{(1)}\ldots m_k^{(k+1)}.a_j,$$
où les $m_j^{(i)}$ sont dans $\M$. Alors, pour tout $j$, 
$$m_j^{(1)}\ldots m_k^{(k+1)}.a_j \in \bigoplus_{l\geq k+1} A(l),$$
donc $x$ ne peut être de degré $k$ : on aboutit à une contradiction, donc $A=(0)$ et donc 
$Q_1,\ldots,Q_n$ engendrent $A_k^G$.\\

{\it Troisième étape.} Montrons que les $Q_i$ sont $(A_1^G)^{\otimes k}$-linéairement indépendants.
Soient $x_1,\ldots,x_n \in (A_1^G)^{\otimes k}$, tels que $x_1Q_1+\ldots+x_nQ_n=0$.
Dans $\frac{A_k}{\M A_k}$, posons, pour tout $j$,
$$Q_j+\M A_k=\sum_{i=1}^m \lambda_{i,j} (P_i+\M A_k).$$
La famille $(Q_j+\M A_k)_{1\leq k\leq n}$ étant libre, la matrice 
$(\lambda_{i,j})_{\substack{1\leq i \leq m\\ 1\leq j \leq n}}\in {\cal M}_{m,n}(\mathbb{C})$ est de rang maximal $n$.
D'autre part, dans $A_k$, posons :
$$Q_j=\sum_{i=1}^m y_{i,j} P_i,$$
où $y_{i,j} \in (A_1^G)^{\otimes k}$ pour tous $i$ et $j$.
 Par unicité des $\lambda_{i,j}$, $y_{i,j}=\lambda_{i,j}+\M A_k$ pour tous $i$ et $j$.
En conséquence, la matrice $(y_{i,j})_{\substack{1\leq i \leq m\\ 1\leq j \leq n}}$ 
est de rang maximal $n$ sur le corps $K=Frac\left((A_1^G)^{\otimes k}\right)$. 
Comme $n \leq m$, son noyau est donc nul. D'autre part,
$$\sum_j x_j Q_j=\sum_{i,j} y_{i,j} x_j P_i,$$
Les $P_i$ étant $(A_1^G)^{\otimes k}$-linéairement indépendants, on obtient pour tout $i$ :
$$\sum_j y_{i,j} x_j=0.$$
D'après ce qui précède, on a donc $x_1=\ldots=x_n=0$. \\

{\it Quatrième étape.} Donc $A_k^G$ est un $(A_1^G)^{\otimes k}$-module libre, de rang fini $n$.
Notons $\breve{Q}_k(h_1,\ldots,h_k)$ la série génératrice de Poincaré-Hilbert des degrés des $Q_i$ ;
il s'agit d'un polynôme, les $Q_i$ étant en nombre fini. 
De plus, la série de Poincaré-Hilbert de $A_k^G$ est :
$$R^G_k(h_1,\ldots,h_k)=\breve{Q}_k(h_1,\ldots,h_k) R^G_1(h_1)\ldots R^G_k(h_k).$$
Donc $\breve{Q}_k(h_1,\ldots,h_k)=Q_k(h_1,\ldots,h_k)$. 
Enfin, le rang $n$ vaut $\breve{Q}_k(1,\ldots,1)$, 
ce qui d'après la proposition \ref{prop3} vaut $|G|^{k-1}$. \\

{\it Dernière étape.} Supposons que $A_k$ soit un $A_k^G$-module libre. Il existe alors un polynôme 
$T_k(h_1,\ldots,h_k)$ tel que $R_k(h_1,\ldots,h_k)=T_k(h_1,\ldots,h_k) R_k^G(h_1,\ldots,h_k)$.
par suite :
$$R_k(h_1,\ldots,h_k)=T_k(h_1,\ldots,h_k) Q_k(h_1,\ldots,h_k) R_1^G(h_1) \ldots R_1^G(h_k).$$
Soient $d_1,\ldots,d_n$ les degrés du groupe de Coxeter $G$ (voir \cite{Bourbaki,Humphreys}). On obtient :
$$\frac{1}{(1-h_1)^n \ldots (1-h_k)^n}
=\frac{T_k(h_1,\ldots,h_k) Q_k(h_1,\ldots,h_k)}
{(1-h_1^{d_1})\ldots (1-h_1^{d_n})\ldots (1-h_k^{d_1})\ldots (1-h_k^{d_n})}.$$
En conséquence, $T_k(h_1,\ldots,h_k)Q_k(h_1,\ldots,h_k)=[d_1]_{h_1}\ldots [d_n]_{h_1}\ldots
[d_1]_{h_k}\ldots [d_n]_{h_k}$.
En considérant la décomposition en polynômes irréductibles de $Q_k(h_1,\ldots,h_k)$, 
on montre que l'on peut écrire $Q_k(h_1,\ldots,h_k)=Q_k^{(1)}(h_1)\ldots Q_k^{(k)}(h_k)$,
où les $Q_k^{(i)}$ sont des polynômes à une variable.
De plus, pour des raisons de symétrie entre les différentes copies de $V$, 
on peut se ramener à $Q_k^{(1)}(h)=\ldots=Q_k^{(k)}(h)=\tilde{Q}_k(h)$.
Comme $Q_k(0,\ldots,0)=1$, on peut supposer que $\tilde{Q}_k(0)=1$.

Comme le rang de $A_k^G$ en tant que $(A_1^G)^{\otimes k}$-module est strictement plus grand que $1$,
 le polynôme $Q_k(h_1,\ldots,h_k)$ n'est pas constant. Soit $\lambda h_1^{\alpha_1}\ldots h_k^{\alpha_k}$
un monôme non constant de ce polynôme, choisi de degré total minimal. A ce monôme correspond $\lambda$
éléments de la famille des générateurs $(Q_1,\ldots,Q_n)$. Si un seul des $\alpha_i$ est non nul,
alors ces générateurs sont dans l'une des copies de $S(V)$ et $G$-invariants, donc dans $(A_1^G)^{\otimes k}$,
ce qui est impossible. Quitte à changer l'indexations des différentes copies de $V$, on peut supposer
que $1\leq \alpha_1 \leq \alpha_2$. Alors $\tilde{Q}_k(h)$ est nécessairement de la forme 
$\tilde{Q}_k(h)=1+\mu h^{\alpha_1}+\ldots$,
où $\mu$ est un scalaire non nul. En développant $Q_k(h_1,\ldots,h_k)$, 
ce dernier polynôme contient le monôme non nul $\mu h_1^{\alpha_1}$, 
ce qui contredit la minimalité du degré de $\lambda h_1^{\alpha_1}\ldots h_k^{\alpha_k}$.
Donc $A_k$ n'est pas libre sur $A_k^G$. $\Box$ \\

{\bf Remarques.} \begin{enumerate}
\item Les trois premières étapes de cette preuve montrent que, si $G$ est un groupe de Coxeter,
$A_k^G$ est un $(A_1^G)^{\otimes k}$-module libre de rang fini, de rang inférieur à $|G|^k$.
\item La dernière étape montre que, si $G$ est un groupe de Coxeter, tel que 
$A_k^G$ soit un $(A_1^G)^{\otimes k}$-module libre de rang $\geq 2$,
alors $A_k$ n'est pas un $A_k^G$-module libre.
\end{enumerate}

\section{Serie de Poincaré-Hilbert des invariants}

\label{sect3}

\subsection{Orthogonal d'un sous-groupe de $(\mathbb{Z}/N\mathbb{Z})^n$}

\begin{defi}
\textnormal{
Soit $K$ un sous-groupe de $(\mathbb{Z}/N\mathbb{Z})^n$. On pose :
$$K^\perp=\left\{(\overline{k}_1,\ldots \overline{k}_n)\in (\mathbb{Z}/N\mathbb{Z})^n\:/\:
\forall (\overline{l}_1,\ldots \overline{l}_n)\in H,\:
\overline{k}_1\overline{l}_1+\ldots +\overline{k}_n\overline{l}_n=\overline{0}\right\}.$$
Il s'agit d'un sous-groupe de $(\mathbb{Z}/N\mathbb{Z})^n$.
De plus, si $H$ est stable sous l'action de $S_n$ par permutation des coordonnées, 
il en est de même pour $H^\perp$.}
\end{defi}

Fixons $k\in \mathbb{N}^*$. Alors $V^{\oplus k}$ a pour base 
$(x_{i,j})_{\substack{1\leq i\leq k \\ 1\leq j\leq n}}$ et l'action
de $G$ sur $V^{\oplus k}$ est donnée par :
$$(\overline{k}_1,\ldots,\overline{k}_n).x_{i,j}=\xi^{k_j} x_{i,j}, \hspace{1cm} \sigma.x_{i,j}=x_{i,\sigma(j)}.$$
Comme $H$ est un sous-groupe distingué de $G$ et que $G/H \approx S_n$, $A_k^G=(A_k^H)^{S_n}$.
Considérons donc d'abord $A_k^H$.
Chaque monôme de $A_k$ engendre un sous-$H$-module de dimension $1$, en conséquence,
$A_k^H=Vect\left(x_{1,1}^{\alpha_{1,1}}\ldots x_{k,n}^{\alpha_{k,n}}\:/\:
x_{1,1}^{\alpha_{1,1}}\ldots x_{k,n}^{\alpha_{k,n}}\mbox{ invariant sous }H \right)$.
De plus, 
\begin{eqnarray*}
&&x_{1,1}^{\alpha_{1,1}}\ldots x_{k,n}^{\alpha_{k,n}}\mbox{ invariant sous }H\\
&\Longleftrightarrow& \forall (\overline{k}_1,\ldots, \overline{k}_n) \in H,\:
\xi^{k_1(\alpha_{1,1}+\ldots+\alpha_{k,1})+\ldots +k_n(\alpha_{1,n}+\ldots+\alpha_{k,n})}=1\\
&\Longleftrightarrow& \forall (\overline{k}_1,\ldots, \overline{k}_n) \in H,\:
\overline{k_1(\alpha_{1,1}+\ldots+\alpha_{k,1})+\ldots +k_n(\alpha_{1,n}+\ldots+\alpha_{k,n})}=\overline{0}\\
&\Longleftrightarrow&(\overline{\alpha_{1,1}+\ldots+\alpha_{k,1}},\ldots,\overline{\alpha_{1,n}+\ldots+\alpha_{k,n}})
\in H^\perp.
\end{eqnarray*}
En conséquence :
\begin{lemme}
\label{lemme6}
$$A_k^H=Vect\left(x_{1,1}^{\alpha_{1,1}}\ldots x_{k,n}^{\alpha_{k,n}}\:/\:
(\overline{\alpha_{1,1}+\ldots+\alpha_{k,1}},\ldots,\overline{\alpha_{1,n}+\ldots+\alpha_{k,n}})
\in H^\perp \right).$$
\end{lemme}

La fin de ce paragraphe est consacrée à la preuve du résultat suivant :
\begin{lemme}
\label{lemme7}
Soit $H$ un sous-groupe de $(\mathbb{Z}/N\mathbb{Z})^n$. Alors :
$$\frac{(\mathbb{Z}/N\mathbb{Z})^n}{H^\perp}\approx
Hom_{\mathbb{Z}}(H,\mathbb{Z}/N\mathbb{Z})\approx H.$$
En conséquence, $|H||H^\perp|=N^n$.
\end{lemme}

{\bf Preuve.} Soit $(e_i)_{1\leq i\leq n}$ la $\mathbb{Z}/N\mathbb{Z}$-base canonique de $(\mathbb{Z}/N\mathbb{Z})^n$.
Il existe une seconde base $(f_i)_{1 \leq i \leq n}$ de $(\mathbb{Z}/N\mathbb{Z})^n$ et des entiers 
$d_1,\ldots,d_n$ tels que :
\begin{enumerate}
\item $d_1\mid d_2\mid \ldots \mid d_n \mid N$ ;
\item $H$ est engendré par $d_1f_1,\ldots, d_nf_n$.
\end{enumerate}

{\it Première étape}. Montrons que l'application suivante est surjective :
$$\rho : \left\{ \begin{array}{rcl}
Hom_{\mathbb{Z}}((\mathbb{Z}/N\mathbb{Z})^n,\mathbb{Z}/N\mathbb{Z}) &\longrightarrow&
Hom_{\mathbb{Z}}(H,\mathbb{Z}/N\mathbb{Z})\\
\phi &\longrightarrow & \phi_{\mid H}.
\end{array}
\right.$$
Soit $\psi \in Hom_{\mathbb{Z}}(H,\mathbb{Z}/N\mathbb{Z})$. Posons $\overline{k}_i=\psi(d_if_i)$ pour tout $i$. 
Comme l'ordre de $d_if_i$ est $N/d_i$, $\overline{k}_i$ est d'ordre divisant $N/d_i$, donc est dans 
$d_i\mathbb{Z}/N\mathbb{Z}$ : posons donc $\overline{k}_i=\overline{d_i l_i}$. 
Alors $\psi$ est la restriction de $H$ du morphisme $\phi$ défini par
$\phi(f_i)=\overline{l_i}$.\\

{\it Deuxième étape}. Considérons l'application suivante :
$$\vartheta  : \left\{ \begin{array}{rcl}
(\mathbb{Z}/N\mathbb{Z})^n&\longrightarrow & Hom_{\mathbb{Z}}(H,\mathbb{Z}/N\mathbb{Z})\\
(\overline{k}_1,\ldots,\overline{k}_n) &\longrightarrow &
 \left\{ \begin{array}{rcl}
H &\longrightarrow & \mathbb{Z}/N\mathbb{Z} \\
(\overline{l}_1,\ldots, \overline{l}_n)&\longrightarrow &
\overline{k}_1\overline{l}_1+\ldots+\overline{k}_n\overline{l}_n.
\end{array}
\right.
\end{array}
\right.$$
Par définition, son noyau est $H^\perp$. Montrons que $\vartheta$ est surjective.
On considère l'application suivante :
$$\theta  : \left\{ \begin{array}{rcl}
(\mathbb{Z}/N\mathbb{Z})^n&\longrightarrow & Hom_{\mathbb{Z}}((\mathbb{Z}/N\mathbb{Z})^n,\mathbb{Z}/N\mathbb{Z})\\
(\overline{k}_1,\ldots,\overline{k}_n) &\longrightarrow &
 \left\{ \begin{array}{rcl}
(\mathbb{Z}/N\mathbb{Z})^n &\longrightarrow & \mathbb{Z}/N\mathbb{Z} \\
(\overline{l}_1,\ldots, \overline{l}_n)&\longrightarrow &
\overline{k}_1\overline{l}_1+\ldots+\overline{k}_n\overline{l}_n.
\end{array}
\right.
\end{array}
\right.$$
Alors $\vartheta=\rho \circ \theta$. D'après la première étape, il suffit de montrer
que $\theta$ est surjectif. On remarque aisément que $(\theta(e_i))_{1\leq i\leq n}$
est une $\mathbb{Z}/N\mathbb{Z}$-base de $Hom_{\mathbb{Z}}((\mathbb{Z}/N\mathbb{Z})^n,\mathbb{Z}/N\mathbb{Z})$,
donc $\theta$ est bijectif. Par suite, 
$\displaystyle \frac{\mathbb{Z}/N\mathbb{Z}}{H^\perp} \approx Hom_{\mathbb{Z}}(H,\mathbb{Z}/N\mathbb{Z})$.\\

{\it Dernière étape.} Montrons que $Hom_{\mathbb{Z}}(H,\mathbb{Z}/N\mathbb{Z}) \approx H$.
Posons $N=d_i d'_i$, pour tout $1\leq i \leq N$. Tout d'abord,
$H=(d_1f_1)\oplus \ldots \oplus (d_nf_n) 
\approx \mathbb{Z}/d'_1\mathbb{Z}\oplus \ldots \oplus \mathbb{Z}/d'_n\mathbb{Z}$.
En conséquence :
$$Hom_{\mathbb{Z}}(H,\mathbb{Z}/N\mathbb{Z})\approx 
Hom_{\mathbb{Z}}\left(\mathbb{Z}/d'_1\mathbb{Z},\:\mathbb{Z}/N\mathbb{Z}\right)\oplus \ldots \oplus Hom_{\mathbb{Z}}\left(d'_n\mathbb{Z},\:\mathbb{Z}/N\mathbb{Z}\right).$$
Il suffit donc de montrer que pour tout $k$ divisant $N$,
$Hom_{\mathbb{Z}}(\mathbb{Z}/k\mathbb{Z},\mathbb{Z}/N\mathbb{Z})$ est cyclique d'ordre $k$.
On considère le morphisme suivant :
$$\phi : \left\{ \begin{array}{rcl}
\mathbb{Z}/k\mathbb{Z}&\longrightarrow & \mathbb{Z}/N\mathbb{Z}\\
\overline{l} &\longrightarrow & \overline{\frac{N}{k}l}.
\end{array}
\right.$$
Ce morphisme est bien défini car $\overline{\frac{N}{k}}$ est d'ordre $k$ dans $\mathbb{Z}/N\mathbb{Z}$.
De plus, $\phi$ est d'ordre $k$ dans $Hom_{\mathbb{Z}}(\mathbb{Z}/k\mathbb{Z},\mathbb{Z}/N\mathbb{Z})$. 
De plus, pour tout morphisme $\psi:\mathbb{Z}/k\mathbb{Z}\longrightarrow \mathbb{Z}/N\mathbb{Z}$,
$\psi(\overline{1})$ est d'ordre divisant $k$, donc de la forme $\overline{\frac{N}{k}l}$ :
par suite, $\psi=l\phi$. Donc $\phi$ engendre $Hom_{\mathbb{Z}}(\mathbb{Z}/k\mathbb{Z},\mathbb{Z}/N\mathbb{Z})$. $\Box$

\subsection{Série formelle de $A_k^G$}

{\bf Notations.}
\begin{enumerate}
\item Soit $n \in \mathbb{N}$. Les partitions de $n$ seront notées
sous la forme $\underline{\alpha}=(\alpha_1,\ldots, \alpha_l)$, avec :
\begin{enumerate}
\item $1\leq \alpha_1\leq \ldots \leq \alpha_l$ ;
\item $\alpha_1+\ldots+\alpha_l=n$.
\end{enumerate}
\item Soit $\underline{\alpha}$ une partition de $n$. Alors $\theta_i(\underline{\alpha})$
est le nombre de $\alpha_j$ égaux à $i$, pour tout $1\leq i \leq n$.\\
\end{enumerate}

{\bf Remarque.} Notons que $l$ dépend de $\underline{\alpha}$. Cependant, pour ne pas alourdir les notations,
nous continuerons à noter $l$ plutôt que $l_{\underline{\alpha}}$.\\

Soit $\sigma\in S_n$ et soit $\underline{\alpha}$ son type. 
Soient $\omega_1,\ldots,\omega_l$ les $\sigma$-orbites, indexées de sorte que 
pour tout $j$, $\omega_j$ soit de cardinal $\alpha_j$. Posons :
$$\chi_k(\sigma)=\sum_{i_1,\ldots,i_k}Tr\left(\sigma_{\mid A^H_k(i_1,\ldots,i_k)}\right)h_1^{i_1}\ldots h_k^{i_k}.$$
Remarquons que $\sigma$ agit par permutation sur les monômes de $A_k^H$. En conséquence,
$\chi_k(\sigma)$ est la série formelle des monômes de $A_k^H$ fixés par $\sigma$.
Pour tout $i \in \{1,\ldots, k\}$, tout $j \in \{1,\ldots,l\}$, posons :
$$x_{i, \omega_j}=\prod_{k \in \omega_j} x_{i,k}.$$
Il s'agit d'un élément de $A_k$ de degré $(0,\ldots,\alpha_j,\ldots,0)$. L'ensemble des monômes de $A_k$
fixés par $\sigma$ est alors :
$$M_k=\left\{x_{1,\omega_1}^{\alpha_{1,1}}\ldots x_{k,\omega_l}^{\alpha_{k,l}}\:/\:
\forall 1\leq i \leq k, \: \forall 1\leq j \leq l, \: \alpha_{i,j} \in \mathbb{N} \right\}.$$
D'après le lemme \ref{lemme6} et par invariance de $H^\perp$ sous l'action de $S_n$, 
un tel monôme est dans $A_k^H$ si, et seulement si,
$$(\underbrace{\alpha_{1,1}+\ldots+\alpha_{k,1}}_{\alpha_1},\ldots,
\underbrace{\alpha_{1,l}+\ldots+\alpha_{k,l}}_{\alpha_l}) \in H^\perp.$$
Le sous-espace engendré par les monômes de $A_k^H$ fixés par $\sigma$ est une sous-algèbre notée $(A_k^H)^\sigma$,
dont la série formelle est $\chi_k(\sigma)$.

\begin{defi}
\textnormal{Soit $\underline{\alpha}$ une partition de $n$. On pose :
\begin{enumerate}
\item $H^\perp_{\underline{\alpha}}=\left\{(\overline{k}_1,\ldots, \overline{k}_l) \in (\mathbb{Z}/N\mathbb{Z})^l \:/\:
(\underbrace{\overline{k}_1,\ldots,\overline{k}_1}_{\alpha_1},\ldots,
\underbrace{\overline{k}_l,\ldots,\overline{k}_l}_{\alpha_l}) \in H^\perp\right\}$.
\item $I_{\underline{\alpha}}(k)=\left\{(\alpha_{i,j})_{\substack{1\leq i \leq k \\1 \leq j\leq l}} \:/\:
0\leq \alpha_{i,j} \leq N-1,\: (\overline{\alpha_{1,1}+\ldots+\alpha_{k,1}}, \ldots, 
\overline{\alpha_{1,l}+\ldots+\alpha_{k,l}})\in H^\perp_{\underline{\alpha}}\right\}.$
\item $\displaystyle P_{\underline{\alpha}}(h_1,\ldots,h_k)
=\sum_{(\alpha_{i,j})\in I_{\underline{\alpha}}(k)}
h_1^{\alpha_1\alpha_{1,1}+\ldots+\alpha_l \alpha_{1,l}}\ldots 
h_k^{\alpha_1\alpha_{k,1}+\ldots+\alpha_l \alpha_{k,l}}.$
\end{enumerate}}
\end{defi}
En particulier, $H_{(1,\ldots,1)}^\perp=H^\perp$.
On remarque en outre que $I_{\underline{\alpha}}(1)$ est en bijection avec $H_{\underline{\alpha}}^\perp$ ;
plus généralement, $|I_{\underline{\alpha}}(k)|=|H_{\underline{\alpha}}^\perp| N^{(k-1)l}$
(on choisit arbitrairement $\alpha_{1,1},\ldots,\alpha_{k-1,1},\ldots,\alpha_{1,l},\ldots,\alpha_{k-1,l}$
et $(\alpha_{k,1},\ldots,\alpha_{k,l})$ est déterminé par l'appartenance à $H_{\underline{\alpha}}^\perp$).\\

Alors, en effectuant la division euclidienne des $\alpha_{i,j}$ par $N$ :
$$(A_k^H)^\sigma=\bigoplus_{(\alpha_{i,j}) \in I_{\underline{\alpha}}(k)}
\mathbb{C}[x_{1,\omega_1}^N,\ldots, x_{k,\omega_l}^N] 
x_{1,\omega_1}^{\alpha_{1,1}} \ldots x_{k,\omega_l}^{\alpha_{k,l}}.$$
En prenant la série formelle de cette sous-algèbre de $A_k$ :
$$\chi_k(\sigma)=\frac{P_{\underline{\alpha}}(h_1,\ldots,h_k)}
{\displaystyle \prod_{i=1}^k \prod_{j=1}^l \left(1-h_i^{N\alpha_j}\right)}.$$
Enfin, avec la proposition \ref{prop2}, comme il existe 
$\displaystyle\frac{n!}{1^{\theta_1(\underline{\alpha})}\ldots n^{\theta_n(\underline{\alpha})}
\theta_1(\underline{\alpha})! \ldots \theta_n(\underline{\alpha})!}$ permutations de type $\underline{\alpha}$ :

\begin{theo}
\label{theo9}
La série de Poincaré-Hilbert de $A_k^G$ vérifie :
$$R^G_k(h_1,\ldots,h_k)=
\sum_{ \mbox{\scriptsize $\underline{\alpha}$ partition de $n$}}
\frac{P_{\underline{\alpha}}(h_1,\ldots,h_k)}{1^{\theta_1(\underline{\alpha})}\ldots n^{\theta_n(\underline{\alpha})}
\theta_1(\underline{\alpha})! \ldots \theta_n(\underline{\alpha})!}
\prod_{i=1}^k \prod_{j=1}^l \frac{1}{1-h_i^{N\alpha_j}}.$$
\end{theo}

On en déduit le corollaire suivant :
\begin{cor}
La série de Poincaré-Hilbert de $A_k^G$ vérifie :
$$\lim_{(h_1,\ldots,h_k) \longrightarrow (1,\ldots,1)}(1-h_1)^n\ldots (1-h_k)^n R^G_k(h_1,\ldots,h_k)
=\frac{1}{|G|}.$$
\end{cor}

{\bf Preuve.} Alors :
\begin{eqnarray*}
&&(1-h_1)^n\ldots (1-h_k)^n R^G_k(h_1,\ldots,h_k)\\
&=& \sum_{\mbox{\scriptsize $\underline{\alpha}$ partition de $n$}}
\frac{P_{\underline{\alpha}}(h_1,\ldots,h_k)}{1^{\theta_1(\underline{\alpha})}\ldots n^{\theta_n(\underline{\alpha})}
\theta_1(\underline{\alpha})! \ldots \theta_n(\underline{\alpha})!}
(1-h_1)^{n-l} \ldots (1-h_k)^{n-l}
\prod_{i=1}^k \prod_{j=1}^l \frac{1}{[N\alpha_j]_{h_i}}.
\end{eqnarray*}
En conséquence, si $\underline{\alpha}\neq(1,\ldots,1)$, alors $l\neq n$ et le terme de la somme
correspondant à $\underline{\alpha}$ tend vers $0$. Par suite :
\begin{eqnarray*}
\lim_{(h_1,\ldots,h_k) \longrightarrow (1,\ldots,1)}(1-h_1)^n\ldots (1-h_k)^n R^G_k(h_1,\ldots,h_k)
&=&\frac{P_{(1,\ldots,1)}(1,\ldots,1)}{n!} \prod_{i=1}^k \prod_{j=1}^n \frac{1}{N}\\
&=&\frac{|I_{(1,\ldots,1)}(k)|}{n!N^{kn}}\\
&=&\frac{|H_{(1,\ldots,1)}^\perp|N^{(k-1)n}}{n!N^{kn}}\\
&=&\frac{|H^\perp|}{n!N^n}.
\end{eqnarray*}
On conclut avec le lemme \ref{lemme7}. $\Box$\\

{\bf Preuve de la proposition \ref{prop3}.} D'après le corollaire précédent,
\begin{eqnarray*}
&&\lim_{(h_1,\ldots,h_k) \longrightarrow (1,\ldots,1)} \frac{R^G_k(h_1,\ldots,h_k)}{R^G_1(h_1)\ldots R^G_1(h_k)}\\
&=&\lim_{(h_1,\ldots,h_k) \longrightarrow (1,\ldots,1)}\frac{(1-h_1)^n \ldots (1-h_k)^n R^G_k(h_1,\ldots,h_k)}
{(1-h_1)^n R^G_1(h_1)\ldots (1-h_k)^n R^G_1(h_k)}\\
&=&\frac{\frac{1}{|G|}}{\frac{1}{|G|} \ldots \frac{1}{|G|}}\\
&=&|G|^{k-1}. \:\Box
\end{eqnarray*}

\section{Exemples des séries infinies de groupes de Coxeter}

\subsection{Groupes symétriques $A_{n-1}$}

On prend $N=2$ et $H=\{(\overline{0},\ldots,\overline{0})\}$.
On obtient ainsi $G=S_n=A_{n-1}$. La représentation associée n'est pas la représentation standard
de $A_{n-1}$, mais la somme directe de la représentation standard $W$ et d'une représentation triviale $T$
de dimension $1$, engendrée par $x_1+\ldots+x_n$. En conséquence, pour tout $k$,
$A_k^G=\mathbb{C}[T^{\oplus k}] \otimes \mathbb{C}[W^{\oplus k}]^G$.
Par suite, le corollaire \ref{theo4} est également vrai pour $A_{n-1}$ agissant sur sa représentation standard.\\

Précisons un peu le théorème \ref{theo9}. On obtient :
\begin{description}
\item[\textnormal{a)}] $H^\perp=(\mathbb{Z}/2\mathbb{Z})^n$.
\item[\textnormal{b)}] Pour toute partition $\underline{\alpha}$, 
$H_{\underline{\alpha}}^\perp=(\mathbb{Z}/2\mathbb{Z})^l$.
\item[\textnormal{c)}] Pour toute partition $\underline{\alpha}$, $I_{\underline{\alpha}}
=\{(\alpha_{i,j})_{\substack{1\leq i\leq k \\ 1\leq j\leq l}}\:/\: \forall i,j, \: 0\leq \alpha_{i,j} \leq 1\}$.
\item[\textnormal{d)}] Pour toute partition $\underline{\alpha}$, 
$$P_{\underline{\alpha}}(h_1,\ldots,h_k)=\sum_{\substack{1\leq i\leq k \\1\leq j \leq l}}
\sum_{0\leq \alpha_{i,j}\leq 1} h_1^{\alpha_1\alpha_{1,1}+\ldots+\alpha_l \alpha_{1,l}}
\ldots h_k^{\alpha_1\alpha_{k,1}+\ldots+\alpha_l \alpha_{k,l}}
=\prod_{\substack{1\leq i\leq k \\1\leq j \leq l}}\left(1+h_i^{\alpha_j}\right).$$
\end{description}
Par exemple :
\begin{enumerate}
\item Pour $n=2$ et $k=2$, $Q_2(h_1,h_2)=h_1 h_2+1$.
\item Pour $n=2$ et $k=3$, $Q_3(h_1,h_2,h_3)=h_1 h_2+h_1 h_3+h_2 h_3 + 1$.
\item Pour $n=3$ et $k=2$, $Q_2(h_1,h_2)=h_1^3 h_2^3+h_1^2 h_2^2+h_1^2 h_2+h_1 h_2^2+h_1 h_2+1$.
\end{enumerate}

\subsection{Groupes de réflexions signées $B_n$}

On prend $N=2$ et $H=(\mathbb{Z}/2\mathbb{Z})^n$.
On obtient ainsi $G=B_n$. La représentation $V$ est la représentation standard de $B_n$.\\

Précisons un peu le théorème \ref{theo9}. On obtient :
\begin{description}
\item[\textnormal{a)}] $H^\perp=\{(\overline{0},\ldots,\overline{0})\}$.
\item[\textnormal{b)}]  Pour toute partition $\underline{\alpha}$, $H_{\underline{\alpha}}^\perp=\{(\overline{0},\ldots,\overline{0})\}$.
\item[\textnormal{c)}] Pour toute partition $\underline{\alpha}$, 
$$I_{\underline{\alpha}}
=\left\{(\alpha_{i,j})_{\substack{1\leq i\leq k \\ 1\leq j\leq l}}\:/\: 
\begin{array}{l}
\forall i,j, \: 0\leq \alpha_{i,j} \leq 1,\\
\alpha_{1,1}+\ldots+\alpha_{k,1},\ldots,\alpha_{1,l}+\ldots+\alpha_{k,l} \mbox{ tous pairs}
\end{array} \right\}.$$
\item[\textnormal{d)}]  Pour toute partition $\underline{\alpha}$, 
$$P_{\underline{\alpha}}(h_1,\ldots,h_k)
=\prod_{j=1}^l \frac{\displaystyle \prod_{i=1}^k
\left(1+h_i^{\alpha_j}\right)+ \prod_{i=1}^k\left(1-h_i^{\alpha_j}\right)}{2}.$$
\end{description}
Par exemple, pour $n=2$ et $k=2$,
$Q_2(h_1,h_2)=h_1^4 h_2^4+h_1^3 h_2^3+h_1^3 h_2+2 h_1^2h_2^2+h_1 h_2^3+h_1 h_2+1$.

\subsection{Groupes de Coxeter $D_n$}

On prend $N=2$ et $H=\{(\overline{k}_1,\ldots,\overline{k}_n) \in (\mathbb{Z}/2\mathbb{Z})^n\:/\:
\overline{k}_1+\ldots+\overline{k}_n=\overline{0}\}$.
On obtient ainsi $G=D_n$. La représentation $V$ est la représentation standard de $D_n$.\\

Précisons un peu le théorème \ref{theo9}. On obtient :
\begin{description}
\item[\textnormal{a)}] $H^\perp=\{(\overline{0},\ldots,\overline{0}), (\overline{1},\ldots,\overline{1})\}$.
\item[\textnormal{b)}]  Pour toute partition $\underline{\alpha}$, 
$H_{\underline{\alpha}}^\perp=\{(\overline{0},\ldots,\overline{0}), (\overline{1},\ldots,\overline{1})\}$.
\item[\textnormal{c)}]  Pour toute partition $\underline{\alpha}$, 
$$I_{\underline{\alpha}}
=\left\{(\alpha_{i,j})_{\substack{1\leq i\leq k \\ 1\leq j\leq l}}\:/\: 
\begin{array}{l}
\forall i,j, \: 0\leq \alpha_{i,j} \leq 1,\\
\alpha_{1,1}+\ldots+\alpha_{k,1},\ldots,\alpha_{1,l}+\ldots+\alpha_{k,l} \mbox{ ont même parité}
\end{array} \right\}.$$
\item[\textnormal{d)}]  Pour toute partition $\underline{\alpha}$, 
$$P_{\underline{\alpha}}(h_1,\ldots,h_k)
=\prod_{j=1}^l \frac{\displaystyle \prod_{i=1}^k
\left(1+h_i^{\alpha_j}\right)+\prod_{i=1}^k\left(1-h_i^{\alpha_j}\right)}{2}
+\prod_{j=1}^l \frac{\displaystyle \prod_{i=1}^k
\left(1+h_i^{\alpha_j}\right)-\prod_{i=1}^k\left(1-h_i^{\alpha_j}\right)}{2}.$$
\end{description}

\subsection{Groupes diédraux $I_2(N)$}

On prend $n=2$, $N$ quelconque et $H=\{(\overline{k}_1,\overline{k}_2) \in (\mathbb{Z}/N\mathbb{Z})^2\:/\:
\overline{k}_1+\overline{k}_2=\overline{0}\}$.
On obtient ainsi le groupe diédral $G=I_N$, d'ordre $2N$. 
La représentation $V$ est la représentation standard de $I_N$.\\

Précisons un peu le théorème \ref{theo9}. On obtient :
\begin{description}
\item[\textnormal{a)}]  $H^\perp=\langle (\overline{1},\overline{1})\rangle$.
\item[\textnormal{b)}]  $H_{(2)}^\perp=\mathbb{Z}/N\mathbb{Z}$ et $H_{(1,1)}^\perp=H^\perp$.
\item[\textnormal{c)}]  On a :
\begin{eqnarray*}
I_{(2)}(k)&=&\{0,\ldots,N-1\}^k,\\
I_{(1,1)}(k)&=&\left\{(\alpha_{i,j})_{\substack{1\leq i \leq k \\1\leq j\leq 2}}\:/\:
\begin{array}{l}
\forall i,j, \: 0\leq \alpha_{i,j} \leq N-1,\\
\alpha_{1,1}+\ldots+\alpha_{k,1}\equiv\alpha_{1,2}+\ldots+\alpha_{k,2}[N]
\end{array} \right\}.
\end{eqnarray*}
\item[\textnormal{d)}]  En conséquence :
\begin{eqnarray*}
P_{(2)}(h_1,\ldots,h_k)&=&\prod_{i=1}^k \left(1+h_i^2+\ldots+h_i^{2(N-1)} \right),\\
P_{(1,1)}(h_1,\ldots,h_k)&=&\sum_{a=0}^{N-1} \left( 
\sum_{b=0}^{N-1} \xi^{ab} 
\frac{\displaystyle \prod_{i=1}^k \left(1+\xi^{b}h_i
+\ldots+ \xi^{b(N-1)}h_i^{N-1}\right)}{N} \right)^2,
\end{eqnarray*}
où $\xi$ est une racine $N$-ième primitive de l'unité.
\end{description}
Par exemple, pour $N=4$ et $k=2$ :
$$Q_2(h_1,h_2)=h_1^4 h_2^4+h_1^3 h_2^3+h_1^3 h_2+2 h_1^2 h_2^2+h_1 h_2^3+h_1 h_2+1.$$

\subsection{Groupes de réflexions complexes $G(de,e,n)$}

Voir par exemple \cite{Geck} pour une description de ces groupes.
On prend ici $N=de$, avec $d$ et $e \in \mathbb{N}^*$ et :
\begin{eqnarray*}
H&=&\{(\overline{k}_1,\ldots,\overline{k}_n) \in (\mathbb{Z}/N\mathbb{Z})^n\:/\:
d(k_1+\ldots+k_n)\equiv 0[N]\}\\
&=&\{(\overline{k}_1,\ldots,\overline{k}_n) \in (\mathbb{Z}/N\mathbb{Z})^n\:/\:
k_1+\ldots+k_n\equiv 0[e]\}.
\end{eqnarray*}
Par suite, $|H|=N^{n-1}d$ et donc $|H^\perp|=e$ par le lemme \ref{lemme7}.
Précisons un peu le théorème \ref{theo9}. On obtient :
\begin{description}
\item[\textnormal{a)}]  $H^\perp=\langle (\overline{d},\ldots,\overline{d}) \rangle$.
\item[\textnormal{b)}] Pour toute partition $\underline{\alpha}$, $H_{\underline{\alpha}}^\perp=
\langle(\overline{d},\ldots,\overline{d})\rangle \subseteq (\mathbb{Z}/N\mathbb{Z})^l$.
\item[\textnormal{c)}] Pour toute partition $\underline{\alpha}$,
$$I_{\underline{\alpha}} =\left\{(\alpha_{i,j})_{\substack{1\leq i\leq k \\ 1\leq j\leq l}}\:/\: 
\begin{array}{l}
\forall i,j, \: 0\leq \alpha_{i,j} \leq 1,\: \exists \lambda \in \mathbb{Z},\\
\alpha_{1,1}+\ldots+\alpha_{k,1}\equiv \ldots \equiv \alpha_{1,l}+\ldots+\alpha_{k,l}\equiv \lambda d[N]
\end{array} \right\}.$$
\item[\textnormal{d)}] Pour toute partition $\underline{\alpha}$, 
$$P_{\underline{\alpha}}(h_1,\ldots,h_k)=
\sum_{\lambda=0}^{e-1} \prod_{j=1}^l
\frac{\displaystyle \sum_{a=0}^{e-1} \xi^{\lambda a} 
\prod_{i=1}^k \left(1+\xi^{a} h_i^{\alpha_j}\right)}{e},$$
où $\xi$ est une racine $e$-ième primitive de l'unité.
\end{description}

\subsection{Un autre exemple, une représentation de $G_2$}

On prend maintenant $N=2$, $n=3$ et $H=\{(\overline{0},\overline{0},\overline{0}),\:
(\overline{1},\overline{1},\overline{1})\}$. On obtient le produit direct $G=H \times S_3$,
isomorphe au groupe diédral $I_2(6)$, ou encore au groupe de Weyl $G_2$. Cependant, la représentation $V$
n'est pas la représentation standard de $G_2$. On obtient facilement que :
$$H^\perp=\{(\overline{0},\overline{0},\overline{0}),\:(\overline{1},\overline{1},\overline{0}),\:
(\overline{1},\overline{0},\overline{1}),\:(\overline{0},\overline{1},\overline{1})\}.$$
Par suite :
\begin{description}
\item[\textnormal{a)}]  $H_{(3)}^\perp=\{\overline{0}\}$ et :
\begin{eqnarray*}
I_{(3)}(k)&=&\left\{(\alpha_i)_{1\leq i \leq k}\:/\:
\begin{array}{l}
\forall i, \: 0\leq \alpha_i \leq 1,\\
\alpha_1+\ldots+\alpha_k \mbox{ est pair}
\end{array} \right\}\: ;\\
P_{(3)}(h_1,\ldots,h_k)&=&
\frac{\displaystyle \prod_{i=1}^k \left(1+h_i^3\right)+\prod_{i=1}^k \left(1-h_i^3\right)}{2}.
\end{eqnarray*}
\item[\textnormal{b)}]  $H_{(1,2)}^\perp=\{(\overline{0},\overline{0}),\:(\overline{0},\overline{1})\}$ et :
\begin{eqnarray*}
I_{(1,2)}(k)&=&\left\{(\alpha_{i,j})_{\substack{1\leq i \leq k\\ 1\leq j \leq 2}}\:/\:
\begin{array}{l}
\forall i,j, \: 0\leq \alpha_{i,j} \leq 1,\\
\alpha_{1,1}+\ldots+\alpha_{k,1} \mbox{ est pair}
\end{array} \right\}\: ;\\
P_{(1,2)}(h_1,\ldots,h_k)&=&
\frac{\displaystyle \prod_{i=1}^k \left(1+h_i\right)+\prod_{i=1}^k \left(1-h_i\right)}{2}
\prod_{i=1}^k(1+h_i^2).
\end{eqnarray*}
\item[\textnormal{c)}]  $H_{(1,1,1)}^\perp=\{(\overline{0},\overline{0},\overline{0}),\:(\overline{1},\overline{1},\overline{0}),\:
(\overline{1},\overline{0},\overline{1}),\:(\overline{0},\overline{1},\overline{1})\}$ et :
\begin{eqnarray*}
P_{(1,1,1)}(h_1,\ldots,h_k)&=&
\left(\frac{\displaystyle \prod_{i=1}^k \left(1+h_i\right)+\prod_{i=1}^k \left(1-h_i\right)}{2}\right)^3\\
&&+3\left(\frac{\displaystyle \prod_{i=1}^k \left(1+h_i\right)+\prod_{i=1}^k \left(1-h_i\right)}{2}\right)
\left(\frac{\displaystyle \prod_{i=1}^k \left(1+h_i\right)-\prod_{i=1}^k \left(1-h_i\right)}{2}\right)^2.
\end{eqnarray*}
\end{description}

On vérifie directement que, pour $k=2$, $Q_2(h_1,h_2)$ est de la forme :
$$Q_2(h_1,h_2)=\frac{(1+h_1h_2)(h_1^6h_2^6 + \ldots)}{(h_1^4+1)(h_2^4+1)}.$$
Ce n'est pas un polynôme. Donc, dans cet exemple,
$\mathbb{C}[V \oplus V]^G$ n'est pas un module libre sur $\mathbb{C}[V]^G \otimes \mathbb{C}[V]^G$.

\bibliographystyle{amsplain}
\bibliography{biblio}

\end{document}